\documentclass[11pt]{amsart}
\usepackage[latin1]{inputenc}
\usepackage{amssymb}
\usepackage{pdfsync}
\usepackage[english]{babel}

\usepackage[pdftex,pagebackref,colorlinks=true,urlcolor=blue,linkcolor=blue,citecolor=blue]{hyperref}
\usepackage{graphicx}

\usepackage[capitalize]{cleveref}
\usepackage{fullpage}
\usepackage{color}
\usepackage{amsmath}
\usepackage{amsfonts}
\usepackage{mathrsfs}
\usepackage{t1enc , graphicx}
\usepackage{verbatim}
\usepackage{bbm}

\usepackage[left= 1.5 in, right= 1.5 in ,top= 1 in, bottom = 2 in]{geometry}
\newcommand{\R}{\mathbb{R}}

\newcommand{\Z}{\mathbb{Z}}
\newcommand{\N}{\mathbb{N}}
\newcommand{\calP}{{\mathcal P}}

\newtheorem{theorem}{Theorem}[section]
\newtheorem{proposition}[theorem]{Proposition}

\newtheorem{lemma}[theorem]{Lemma}

\theoremstyle{definition}

\theoremstyle{remark}

\newtheorem{remark}[theorem]{Remark}

\author{Anh Le}
\author{Th\'ai Ho\`ang L\^e}

\address{Department of Mathematics\\
	Northwestern University\\
	2033 Sheridan Road, Evanston, IL 60208, USA}   
\email{anhle@math.northwestern.edu  }

\address{Department of Mathematics\\
    University of Mississippi\\
    University, MS 38677, USA}
\email{leth@olemiss.edu}

\title{Multidimensional configurations in the primes with shifted prime steps}

\begin{document}
\begin{abstract}
	Let $\mathcal{P}$ denote the set of primes. For a fixed dimension $d$, Cook-Magyar-Titichetrakun, Tao-Ziegler and Fox-Zhao independently proved that any subset of positive relative density of $\mathcal{P}^d$ contains an arbitrary linear configuration. In this paper, we prove that there exists such configuration with the step being a shifted prime (prime minus $1$ or plus $1$). 
\end{abstract}
\maketitle
\section{Introduction}

\subsection{History and statement of the result}

For a set $A \subset \Z^d$, we define its upper density by $\overline{d}(A) = \limsup_{N \rightarrow \infty} \frac{|A \cap [1,N]^d|}{N^d}$. We say $A$ is dense if $\overline{d}(A)>0$.
The celebrated Szemer\'edi's theorem \cite{Szemeredi75} states that any dense subset of the integers must contain arbitrarily long arithmetic progressions. Later, Furstenberg \cite{Furstenberg77} provided an ergodic theoretic proof of this result. After this proof, many far-reaching generalizations have been obtained. For example, Furstenberg-Katznelson \cite{Furstenberg-Katznelson-1978} proved a multidimensional generalization of Szemer\'edi's theorem, namely that for any vectors $h_1,\ldots, h_t \in \Z^d$, any dense subset of $\Z^d$ contains a configuration of the form $\{a + r h_1, \ldots, a + r h_t \}$ for some $a \in \Z^d$ and $r \in \Z$ nonzero. Bergelson-Leibman \cite{Bergelson_Leibman96} generalized Furstenberg-Katznelson's result to polynomial configurations. That is to say: Given vectors $h_1, \ldots, h_t \in \Z^d$ and polynomials $P_1, \ldots, P_t \in \Z[x]$ without constant terms, any dense subset of $\Z^d$ must contain a configuration of the form $\{ a +P_1(r) h_1, \ldots, a + P_t(r) h_t \}$ for some $a \in \Z^d$ and $r \in \Z$ nonzero.

In \cite{Frantzikinakis-Host-Kra-07} and  \cite{Frantzikinakis-Host-Kra-13}, Frantzikinakis-Host-Kra showed that the ``step'' $r$ in the theorems of Furstenberg-Katznelson and Bergelson-Leibman can be taken to be a shifted prime, i.e. a number of the form $p-1$ where $p$ is a prime\footnote{The same result is true with $p+1$ in place of $p-1$. Simple counterexamples show that these are the only translates of the primes enjoying this property.}, using what we now term the comparison method. These results were also proved independently by Wooley-Ziegler \cite{Wooley_Ziegler12} (in the polynomial case) and Bergelson-Leibman-Ziegler \cite{Bergelson-Leibman-Ziegler-11} (in the multidimensional case).

Let $\calP$ denote the set of primes. For a set $A \subset \calP^d$, we define its relative upper density by $\overline{d}_{\calP^d}(A) = \limsup_{N \rightarrow \infty} \frac{|A \cap [1,N]^d|}{|\calP^d \cap [1,N]^d|}$. We say $A$ is dense in $\calP^d$ if $\overline{d}_{\calP^d}(A)>0$. The Green-Tao theorem \cite{Green_Tao08} states that any dense subset of $\calP$ must contain an arbitrarily long arithmetic progression. Tao-Ziegler \cite{Tao_Ziegler08} generalized this to polynomial configurations, namely any dense subset of $\calP$ must contain a configuration $\{ a +P_1(r), \ldots, a + P_t(r) \}$ for some $a, r \in \Z, r \neq 0$, where $P_1, \ldots, P_t \in \Z[x]$ are given polynomials without constant terms. Regarding multidimensional configurations, Cook-Magyar-Titichetrakun \cite{Cook-Magyar-Titichetrakun-2015} and Tao-Ziegler \cite{Tao-Ziegler-2015} independently proved a prime version of the Furstenberg-Katznelson theorem. That is to say for any $h_1, h_2, \ldots, h_k \in \Z^d$, any dense subset of $\mathcal{P}^d$ must contain a configuration $\{ a + r h_1, a + r h_2, \ldots, a + r h_k \}$ for some $a \in \Z^d$ and $r \in \Z$ nonzero. Shortly after, Fox-Zhao \cite{Fox-Zhao-2015} came up with a very short proof of the same result.

In \cite{Le-Wolf-2014}, Wolf and the second author ``completed the square'' by proving a hybrid of Tao-Ziegler \cite{Tao_Ziegler08} and Wooley-Ziegler \cite{Wooley_Ziegler12}'s results, namely that in polynomial configurations $\{ a +P_1(r), \ldots, a + P_t(r) \}$ in dense subsets of the primes, one can require the step $r$ to be of the form $p-1$ (or $p+1$), where $p$ is prime. Their proof relies on the comparison method of Frantzikinakis-Host-Kra. In this paper, we complete another square by showing that the step $r$ in the multidimensional result of Cook-Magyar-Titichetrakun, Tao-Ziegler and Fox-Zhao can also be taken to be of the form $p-1$ or $p+1$. To be precise, we will prove the following:

\begin{theorem}
    \label{thm:sept-22-1}
    Let $A \subseteq \mathcal{P}^d$ be of positive relative upper density. Then for every $h_1, \ldots, h_k \in \mathbb{Z}^d$, there exists $a \in \mathbb{Z}^d$ and $p \in \mathcal{P}$ such that $a + (p-1) h_1, \ldots, a + (p-1) h_k \in A$. The same is true if we replaced $p-1$ by $p+1$.
\end{theorem}
In fact, we are able to find a lower bound for the number of such configurations. Again, simple examples show that $\mathcal{P} - 1$ and $\mathcal{P} + 1$ are the only translates of the primes having this property.

\subsection{Idea of the proof} In \cite{Le-Wolf-2014}, the authors used Tao-Ziegler's theorem on polynomial configurations in the primes as a black box and employed Frantzikinakis-Host-Kra's comparison method. In the current situation, our implementation of the comparison method is less straightforward due to different kinds of averages being considered. Our proof relies on Tao-Ziegler \cite{Tao-Ziegler-2015}'s weighted Furstenberg correspondence principle, which we now recall. With the setup as in \cref{thm:sept-22-1}, Tao and Ziegler constructed a measure preserving system $(X, \mathcal{B}, \mu, (T_h)_{h \in \Z^d})$ and a set $E \in \mathcal{B}$ such that for $m \in \N$, the expression
\begin{equation}
    \label{eq:nov-14-1}
    \mu(T_{mh_1} E \cap \ldots \cap T_{m h_k} E)
\end{equation}
essentially captures the density of homothetic copies of  $\{h_1, \ldots, h_k\}$ in $A$. The existence of $m$ that makes \eqref{eq:nov-14-1} positive is then guaranteed by the Furstenberg-Katznelson theorem \cite{Furstenberg-Katznelson-1978}. The whole idea is reminiscent of that of Furstenberg \cite{Furstenberg77} when proving Szemer{\'e}di's theorem, in which he constructed a system and a set arising from a dense subset of $\Z$ (instead of $\mathcal{P}^d$).

To show the step of the configurations found in \cite{Tao-Ziegler-2015} can be restricted to shifted primes, one is tempted to use the Franzikinakis-Host-Kra's comparison method. However, this approach does not produce configurations as desired. In contrast to the measure built in \cite{Furstenberg77}, the measure in \cite{Tao-Ziegler-2015} contains an extra average on a parameter $r$ (see \cref{sec:tao-ziegler-summary}). Hence with this approach, we end up with configurations of the form 
\begin{equation}
    \label{eq:nov-14-2}
    a + m r h_1, \ldots, a + m r h_k
\end{equation}
with $m = p - 1$ for some prime $p$ and some $r \in \Z$ which we have no control of. The appearance of $r$ arises from the use of Varnavides' trick \cite{Varnavides}. We remark that Fox-Zhao's proof \cite{Fox-Zhao-2015} also used Varnavides' trick. While their proof is very simple, the number of configurations it provides is too small for our purpose.

Instead, to prove \cref{thm:sept-22-1}, we first use Furstenberg and Katznelson's theorem to fix an $m$ that makes \eqref{eq:nov-14-1} positive. Only after that, the comparison method is implemented, taking advantage of the additional average. As a result, we can obtain configurations as \eqref{eq:nov-14-2} with $r = (p-1)/m$ for some prime $p$. This produces the desired result.  

The detailed proof will be presented in \cref{sec:main-proof} contingent upon Proposition \ref{prop:1}, which is proved in \cref{sec:prop-3.1}. In the next section, we set up the notation and necessary background. 

\subsection{Acknowledgment} The authors would like to thank Bryna Kra for useful feedback. The first author also thanks the hospitality of the University of Mississippi where part of this work was carried out. The second author is supported by National Science Foundation Grant DMS-1702296.

\section{Preliminaries}
\label{sec:tools}
\subsection{Notation}
\label{sec:notation}
For $N \in \mathbb{N}$, let $[N]$ denote the set $\{1, 2, \ldots, N\}$. 

For a finite set $S$ and a function $f: S \to \R$, let $\mathbb{E}_{n \in S} f(n)$ denote the average $\frac{1}{|S|} \sum_{n \in S} f(n)$.

For two functions $f, g: \N \to \R$, by writing $f(n) \ll g(n)$, we mean there is a constant $c > 0$ such that $|f(n)| \leq c g(n)$ for all $n \in \N$.  

Let $o_{n \to \infty}(1)$ or $o(1)$ denote a function that goes to $0$ as $n$ tends to infinity.

Let $\omega_n$ be a function growing slowly with $n$ (to be determined later). Then we define
\begin{equation*}
    W_n = \prod_{p \leq \omega_n, p \in \mathcal{P}} p.
\end{equation*}

For $-W_n < b < W_n$ coprime to $W_n$, define the $W$-tricked von Mangoldt function as
\begin{equation*}
    \Lambda_{W_n, b}(m) = \frac{\phi(W_n)}{W_n} \Lambda'(W_n m + b)
\end{equation*}
where $\Lambda'$ is the modified von Mangoldt function
\begin{equation*}
    \Lambda'(m) = \begin{cases} \log m \, \mbox{ if } m \in \mathcal{P} \\
    0 \, \mbox{ otherwise.}
    \end{cases}
\end{equation*}

\subsection{Furstenberg $\mathbb{Z}^d$-system}
Let $X = 2^{\mathbb{Z}^d} = \{B: B \subseteq \mathbb{Z}^d\}$. Let $\mathcal{B}$ be the $\sigma$-algebra on $X$ which is generated by basic cylinder sets $\{B \in X: b \in B\}$ for some $b \in \mathbb{Z}^d$. This space has an $\Z^d$-action $(T_h)_{h \in \mathbb{Z}^d}$ defined by $T_h B = B + h$ for all $h \in \mathbb{Z}^d$ and $B \in X$. Then $(X, \mathcal{B}, (T_h)_{h \in \mathbb{Z}^d})$ is a topological $\mathbb{Z}^d$-system, which we call the Furstenberg $\mathbb{Z}^d$-system.
\subsection{The Furstenberg-Katznelson theorem}
To provide a uniform lower bound for the number of desired configurations, we will utilize the following uniform version of the Furstenberg-Katznelson theorem by Bergelson, Host, Mccutcheon and Parreau. 
\begin{theorem}[{\cite[Theorem 2.1 (ii)]{Bergelson_Host_McCutcheon_Parreau2000}}] \label{th:fk}
For any $k$ and for any $\delta >0$, there exists a constant $c(\delta, k)>0$ depending only on $\delta$ and $k$ such that the following holds. For any $k$ commuting measure preserving transformations $T_1, \ldots, T_k$ of a probability measure space $(X, \mathcal{B}, \mu)$ and any $A \in \mathcal{B}$ with $\mu(A) = \delta $, there exists a positive integer $n$ such that $\mu(T_1^{-n} A \cap \cdots \cap T_k^{-n} A)> c(\delta,k)$.
\end{theorem}

\subsection{Summary of Tao and Ziegler's construction}
\label{sec:tao-ziegler-summary}
Our proof relies heavily on the proof of Tao and Ziegler \cite{Tao-Ziegler-2015}. Since many detailed constructions in that proof will be used, we summarize them here for the convenience.

From now on we fix a dimension $d$, as well as a dense subset $A \subseteq \mathcal{P}^d$. There exists a $\delta > 0$ and a sequence $N_n$ going to infinity as $n \to \infty$ such that
\[
    |A_n| \geq \delta |\mathcal{P} \cap [N_n]|^d
\]
where $A_n = A \cap [N_n]^d$. By deleting a small number of elements from $A_n$ (and reducing $\delta$ accordingly), we may assume $A_n \in [\delta' N_n, (1- \delta') N_n]^d$ for some $0 < \delta' < 1/2$. Note that in \cite{Tao-Ziegler-2015}, the authors only removed the primes in $[0, \delta' N_n]$. The reason for this removal is to make sure when $n$ is large enough, all primes in $A_n$ are coprime to $W_n$. Then they used the pigeonhole principle to choose $b_1, \ldots, b_d$ satisfying \eqref{eq:oct-12-1} later. Here, for our purpose, we also delete the primes in $[(1- \delta') N_n, N_n]$. It is because later, we will need to shift the variable $a$ in $A_n$. This deletion makes sure the shift does not move $a$ out of $[N_n]$. 

Let $\omega_n \ll \log \log \log N_n$ be a sufficiently slowly growing function which will be chosen later. Define $W_n = \prod_{p \leq \omega, p \in \mathcal{P}} p $ as in \cref{sec:notation}. For $-W_n < b_1, \ldots, b_d < W_n$ coprime to $W_n$, denote $b = (b_1, \ldots, b_d) \in \mathbb{Z}^d$. Define $N'_n = \lfloor N_n/W_n \rfloor$ and
\begin{equation}
    A'_n = \{a \in [N'_n]^d: W_na + b \in A_n\}
\end{equation}

Let $M_n = o(N'_n)$ be a sequence of natural numbers. For each $a \in A'_n$ and $r \in M_n$, let $B_{a, r, n} = \{b \in \mathbb{Z}^d: a + rb \in A'_n\}$.
For each finite set $\Omega \subset \mathbb{Z}^d$ and $n \in \mathbb{N}$, define a measure $\mu_{\Omega, n}$ on the Furstenberg $\mathbb{Z}^d$-system $X$ by
\[
    \mu_{\Omega, n} = \mathbb{E}_{a \in [N'_n]^d} \mathbb{E}_{r \in [M_n]} \delta_{B_{a,r,n}} \prod_{i=1}^d \prod_{c_i \in \Omega_i} \Lambda_{W_n, b_i} (a_i + c_i r)
\]
where $\Omega_i$ is the projection of $\Omega$ on the $i$-th coordinate, $a = (a_1, a_2, \ldots, a_d)$ and $\delta$ is the delta mass, i.e.
\begin{equation}
    \delta_B(S) = \begin{cases} 1 \, \mbox{ if }B \subset S \\
    0 \, \mbox{ otherwise}
    \end{cases}
\end{equation}
for any subsets $B$ and $S$ of $\mathbb{Z}^d$. Let $E$ be the basic cylinder set 
\[
    E = \{B \in X: 0_{\mathbb{Z}^d} \in B\}
\]
Then for every $h_1, \ldots, h_k \in \mathbb{Z}^d$, by definition of the measure $\mu_{\Omega, n}$, one has
\begin{multline}
\label{eq:mar-15-1}
    \mu_{\Omega, n}(T_{h_1} E \cap \cdots \cap T_{h_k} E) = \\
    \mathbb{E}_{a \in [N'_n]^d} \mathbb{E}_{r \in [M_n]} \prod_{j=1}^k 1_{A'_n}(a + rh_j) \prod_{i=1}^d \prod_{c_i \in \Omega_i} \Lambda_{W_n, b_i}(a_i + c_i r)
\end{multline}

In \cite{Tao-Ziegler-2015}, it is shown that we can choose $\omega_n, b_1, \ldots, b_d$ and $M_n$ such that $\Lambda_{W_n, b_i}$ for $1 \leq i \leq d$ satisfy the linear forms condition (see \cref{prop:linear-form-condition} below) and 
\begin{equation}
\label{eq:oct-12-1}
    \lim_{n \to \infty} \mathbb{E}_{a \in [N_n']^d} 1_{A_n'}(a) \prod_{i=1}^d \Lambda_{W_n, b_i}(a_i) \geq \delta.
\end{equation}
Then we define a measure $\mu$ on $X$ by
\begin{equation}
\label{eq:sept-22-2}
    \mu = p-\lim_{l \to \infty} ( p-\lim_{n \to \infty} \mu_{[-l,l]^d, n})
\end{equation}
where $p-\lim$ is a fixed Banach limit functional, i.e. a linear functional extending the standard limit functional $\lim$ on convergent sequences such that
\begin{equation}
    \liminf_{n \to \infty} x_n \leq p-\lim_{n \to \infty} x_n \leq \limsup_{n \to \infty} x_n.
\end{equation}

It is easy to see that the measure $\mu$ is well-defined and is a probability measure on $(X, \mathcal{B})$. The crucial point is that Tao and Ziegler used the linear forms condition to prove $\mu$ is invariant under $T_h$ for all $h \in \mathbb{Z}^d$ and to prove the compatibility property (\cite[Proposition 2.7]{Tao-Ziegler-2015}). Among other things, this compatibility implies that for every $h_1, \ldots, h_k \in \mathbb{Z}^d$ 
\begin{equation}
\label{eq:oct-13-1}
    \mu_{\Omega, n}(T_{h_1} E \cap \ldots \cap T_{h_k} E) = \mu_{\{h_1, \ldots, h_k\}, n} (T_{h_1} E \cap \ldots \cap T_{h_k} E) + o(1)
\end{equation}
as long as $\{h_1, \ldots, h_k\} \subseteq \Omega$.
 
With the measure $\mu$, one gets $(X, \mathcal{B}, \mu, (T_h)_{h \in \mathbb{Z}^d})$ is a measure preserving $\mathbb{Z}^d$-system. From \eqref{eq:oct-12-1} and \eqref{eq:oct-13-1},
\begin{equation}
    \mu(E) = p-\lim_{n \to \infty} \mu_{\{0\}, n}(E) = p-\lim_{n \to \infty} \mathbb{E}_{a \in [N'_n]} 1_{A'_n}(a) \prod_{i=1}^d \Lambda_{W_n, b_i}(a_i) \geq \delta.
\end{equation}

\subsection{Linear forms condition}
As mentioned earlier, the construction of the measure $\mu$ relies on the following property of $W$-tricked von Mangoldt functions.
\begin{proposition}[cf. {\cite[Proposition 2.5]{Tao-Ziegler-2015}}]
    \label{prop:linear-form-condition}
    Let $d \geq 1$ be a fixed dimension. Let $N'_n$ be a sequence of natural numbers going to infinity as $n \to \infty$. Then we can choose sequences $M_n = o(N'_n)$, $H_n = o(M_n)$ and $\omega_n$ such that the following is true:
    
    Let $m, k_1, \ldots, k_d \geq 0$ be natural numbers. For $1 \leq i \leq d, 1 \leq j \leq k_i$, let $\phi_{i,j}: \mathbb{Z}^{d+m} \to \mathbb{Z}$ be linear forms with integer coefficients and pairwise linearly independent. One has
    \begin{equation} \label{eq:linearforms}
        \mathbb{E}_{a \in [N'_n]^d} \mathbb{E}_{r \in \prod_{j=1}^m L_{j,n}} \prod_{i=1}^d \prod_{j=1}^{k_i} \nu_{i,j,n} (\phi_{i,j} (a, r)) = 1 + o(1)
    \end{equation}
    for all subintervals $L_{j,n}$ of $[-M_n, M_n]$ of length greater than $H_n$ and $\nu_{i,j,n} = \Lambda_{W_n, b}$ for some $(b, W_n) = 1$. 
\end{proposition}
\begin{remark}
    The linear forms condition in \cite[Proposition 2.5]{Tao-Ziegler-2015} is weaker than the version presented here in the sense that it is restricted to linear forms of the form $\phi_{i,j}(a, r) = a_i + \psi_{i,j}(r)$ for some linear forms $\psi_{i,j}:\mathbb{Z}^m \to \mathbb{Z}$ where $a = (a_1, \ldots, a_d)$. The generalization here is needed for \cref{prop:1} and its proof is identical to \cite[Page 17]{Tao-Ziegler-2015}. We summarize below for completeness.   
\end{remark}
\begin{proof}
    First note that the collection of $(m, k_1, \ldots, k_d, \phi_{i,j})$ is countable. Fix an arbitrary enumeration of them. By \cite[Theorem 5.1]{Green_Tao10}, for $M \in \mathbb{N}$, there exists $\omega = \omega^{(M)}$ and $n^{(M)}$ such that
    \begin{equation}
        \left| \mathbb{E}_{a \in [N'_n]^d} \mathbb{E}_{r \in \prod_{j=1}^m L_{j,n}} \prod_{i=1}^d \prod_{j=1}^{k_i} \nu_{i,j,n}(\phi_{i,j}(a,r)) - 1 \right| < \frac{1}{M}
    \end{equation}
    for the first $M$ choices of $m, k_1, \ldots, k_d, \phi_{i,j}$ as long as $L_{j,n}$ are subintervals of $[-N'_n/M, N'_n/M]$ of length greater than $N'_n/M^2$ and $n \geq n^{(M)}$. 
    
    We can choose $n^{(M)}$ so that the sequence $n^{(M)}$ is increasing in $M$. For $M \in \mathbb{N}$ and $n^{(M)} \leq n < n^{(M+1)}$, let $\omega_n = \omega^{(M)}$, $M_n = N'_n/M$ and $H_n = N'_n/M^2$. These choices satisfy the conclusion of \cref{prop:linear-form-condition}.
\end{proof}

\subsection{Other tools}
Similar to \cite{Le-Wolf-2014}, we need following elementary lemmas.

\begin{lemma}[Cauchy-Schwarz]
Let $A, B$ be finite sets, $f, F$ be functions on $A$ and $g$ be a function on $A \times B$. If $|f| \leq F$ pointwise, then 
\begin{equation*}
    |\mathbb{E}_{a \in A, b \in B} f(a) g(a,b)|^2 \leq \mathbb{E}_{a \in A} F(a) \times \mathbb{E}_{a \in A} F(a) |\mathbb{E}_{b \in B} g(a,b)|^2
\end{equation*}
\end{lemma}

\begin{lemma}[van der Corput]
Let $(x_n)_{n \in \mathbb{Z}}$ be complex-valued sequence satisfying $x_n = 0$ outside of interval $[N]$. Then
\begin{equation*}
    \label{eq:oct-10-1}
    |\mathbb{E}_{n \in [N]} x_n|^2 \ll \mathbb{E}_{|h|< N} \mathbb{E}_{n \in [N]} x_n x_{n+h}
\end{equation*}
\end{lemma}

\section{Proof of Theorem \ref{thm:sept-22-1}}
\label{sec:main-proof}
In this section, we prove \cref{thm:sept-22-1} conditional on a proposition whose proof will be presented in the next section. 
Let the setup be as in \cref{sec:tao-ziegler-summary}.
Fix $h_1, h_2, \ldots, h_k \in \mathbb{Z}^d$. By \eqref{eq:sept-22-2}, for every $m \in \mathbb{N}$, one has
\begin{equation}
    \mu(T_{m h_1} E \cap \ldots \cap T_{m h_k} E) = p-\lim_{l \to \infty} p-\lim_{n \to \infty} \mu_{[-l, l]^d, n} (T_{m h_1} E \cap \ldots \cap T_{m h_k} E)
\end{equation}
By \eqref{eq:oct-13-1} (compatibility property), the second $p-\lim$ is redundant if $[-l,l]^d$ is replaced by $\{m h_1, \ldots, m h_k\}$. To be precise, we have
\begin{equation}
\label{eq:sept-22-4}
    \mu(T_{m h_1} E \cap \ldots \cap T_{m h_k} E) = p-\lim_{n \to \infty} \mu_{\{m h_1, \ldots, m h_k\}, n} (T_{m h_1} E \cap \ldots \cap T_{m h_k} E)
\end{equation}
Since $(X, \mathcal{B}, \mu, (T_h)_{h \in \mathbb{Z}^d})$ is a measure preserving $\mathbb{Z}^d$-system and $\mu(E) = \delta > 0$, Theorem \ref{th:fk} implies that there exist $c(\delta, k) >0$ and $m \in \mathbb{N}$ such that
\begin{equation}
    \mu(T_{m h_1} E \cap \ldots \cap T_{m h_k} E) > c(\delta,k).
\end{equation}
By \eqref{eq:sept-22-4},
\[
    p-\lim_{n \to \infty} \mu_{{m h_1, \ldots, m h_k},n} (T_{m h_1} E \cap \ldots T_{m h_k} E) > c(\delta, k)
\]
This means
\begin{equation}
\label{eq:sept-22-5}
    \mu_{\{m h_1, \ldots, m h_k\}, n} (T_{m h_1} E \cap \ldots \cap T_{m h_k} E) > c(\delta, k) - o(1)
\end{equation}
By \eqref{eq:mar-15-1}, the inequality \eqref{eq:sept-22-5} implies
\begin{equation}
    \label{eq:2}
    \mathbb{E}_{a \in [N'_n]^d} \mathbb{E}_{r \in [M_n]} \prod_{j=1}^k 1_{A'_n}(a + rmh_j) \prod_{i=1}^d \prod_{c_i \in \Omega_i} \Lambda_{W_n, b_i}(a_i + c_i r) > c(\delta, k) - o(1)
\end{equation}
where $\Omega_i \subset \mathbb{Z}$ is the projection of the set $\Omega =\{m h_1, \ldots, m h_k\}$ on the $i$-th coordinate.

The following proposition is our version of the comparison method.
\begin{proposition}
    \label{prop:1}
    Let everything be as before. Then 
    \begin{multline}
    \label{eq:sept-21-7}
        \mathbb{E}_{r \in [M_n]} (\Lambda_{W_n, 1}(mr) - 1) \mathbb{E}_{a \in [N'_n]^d} \prod_{j=1}^k 1_{A'_n}(a + rmh_j) \prod_{i=1}^d \prod_{c_i \in \Omega_i} \Lambda_{W_n, b_i}(a_i + c_i r) \\= o(1).
    \end{multline}
\end{proposition}

Assume we already had Proposition \ref{prop:1}. Then \eqref{eq:sept-21-7} together with \eqref{eq:2} imply that 
\begin{multline}
\label{eq:nov-12-1}
     \mathbb{E}_{r \in [M_n]} \mathbb{E}_{a \in [N'_n]^d} \Lambda_{W_n, 1}(mr) \prod_{j=1}^k 1_{A'_n}(a + rmh_j) \prod_{i=1}^d \prod_{c_i \in \Omega_i} \Lambda_{W_n, b_i}(a_i + c_i r) > \\ 
     c(\delta, k) - o(1).
\end{multline}
Hence for any sufficiently large $n$, the set $A'_n$ contains a configuration of the form $a + rm h_1, \ldots, a + rm h_k$ with $W_nmr + 1 \in \mathcal{P}$. This is equivalent to saying that $A_n$ contains $W_na + b + W_n mr h_j$ for $1 \leq j \leq k$ with $W_n mr \in \mathcal{P} - 1$. This finishes our proof of Theorem \ref{thm:sept-22-1} for the step being $p-1$ for some $p \in \mathcal{P}$. For the case $p+1$, simply replace $\Lambda_{W_n,1}(mr) - 1$ by $\Lambda_{W_n, -1}(mr) -1$ in \cref{prop:1}. The rest of the proof remains the same. 

In fact, \eqref{eq:nov-12-1} gives us a lower bound for number of the pairs $(\tilde{a}, p) \in [N_n]^d \times [W_n m M_n + 1]$ satisfying $p \in \mathcal{P}$ and $\tilde{a} + (p-1) h_j \in A$ for all $1 \leq j \leq k$. This number is not less than the number of pairs $(a, r) \in [N_n']^d \times [M_n]$ such that $W_n m r + 1\in \mathcal{P}$ and $a + m r h_j \in A_n'$. By \eqref{eq:nov-12-1}, the latter is greater than
\begin{equation}
    \label{eq:nov-12-2}
    (c(\delta, k) - o(1)) \times \frac{M_n N_n'^d}{\log(W_n m M_n + 1) \log^{\sum_{i=1}^d |\Omega_i|} N_n} \times \left(\frac{W_n}{\phi(W_n )}\right)^{1 + \sum_{i=1}^d |\Omega_i|}.
\end{equation}
Since $N_n' = \lfloor N_n/W_n \rfloor$ and $W_n = o(\log(M_n))$, we get \eqref{eq:nov-12-2} is equal to
\begin{equation}
    \label{eq:nov-13-1}
    (c(\delta, k) - o(1)) \times \frac{M_n N_n^d}{\log M_n \log^{\sum_{i=1}^d |\Omega_i|} N_n} \times \frac{W_n^{1 -d + \sum_{i=1}^d |\Omega_i|}}{\phi(W_n )^{1 + \sum_{i=1}^d |\Omega_i|}}.
\end{equation}

\section{Proof of Proposition \ref{prop:1}}
\label{sec:prop-3.1}

The rest of our paper is devoted to proving \cref{prop:1}. The strategy is similar to \cite[Proposition 1]{Le-Wolf-2014} which in turn is inspired by the method in \cite{Frantzikinakis-Host-Kra-07, Frantzikinakis-Host-Kra-13}. The idea here is that after a finite number of applications of Cauchy-Schwarz and van der Corput, the left hand side of \eqref{eq:sept-21-7} is bounded by an expression consisting entirely of $W$-tricked von Mangoldt functions. This expression then tends to zero by the linear forms condition.

\subsection{A toy example}
To illustrate the ideas of the proof, we first work with the following special case. The proof of the general case is not much different. 

In this model case, we take $d = 2$, $h_1 = (1, 1)$ and $h_2 = (2, 1)$. Then the projection of $\Omega = \{h_1, h_2\}$ onto the first coordinate is $\Omega_1 = \{1, 2\}$ and onto the second coordinate is $\Omega_2 = \{1\}$. We need to show:
\begin{multline}
\label{eq:sept-21-1}
    \mathbb{E}_{r \in [M_n]} (\Lambda_{W, 1}(mr) - 1) \mathbb{E}_{(a_1, a_2) \in [N'_n]^2} 1_{A'_n}((a_1 + mr, a_2 + mr)) 1_{A'_n}((a_1 + 2mr, a_2 + mr)) \times \\
    \times \Lambda_{W, b_1}(a_1 + mr) \Lambda_{W, b_1}(a_1 + 2mr) \Lambda_{W, b_2}(a_2 + mr) = o(1).
\end{multline}

From now on, for simplicity, we write $W$ instead of $W_n$. Likewise, we will assume $m = 1$ since doing so significantly simplifies the notation while not affecting the proof. 

For $n \in \mathbb{N}$ and $a = (a_1, a_2) \in \mathbb{Z}^2$, define following functions:
\begin{itemize}
    \item $\nu(n) = \Lambda_{W, 1}(n)$
    
    \item $\nu_{1, 1}(a) = \Lambda_{W, b_1}(a_1)$ and $\nu_{1, 2}(a) = \Lambda_{W, b_2}(a_2)$.
    
    \item $\nu_{2, 1}(a) = \Lambda_{W, b_1}(a_1)$ and $\nu_{2, 2} (a) = 1$.
    
    \item $\theta_1(a) = \nu_{1, 1}(a) \nu_{1, 2}(a)$   
    
    \item $\theta_2(a) = \nu_{2, 1}(a) \nu_{2, 2}(a)$  
\end{itemize}
Note that the way we define $\nu_{i,j}$ $(i, j = 1, 2)$ depends on the set $\{h_1, h_2\}$. Since the second coordinate of $h_2$ coincides with the second coordinate of $h_1$, we define $\nu_{2,2}(a) = 1$. On the other hand, if they did not coincide, we would define $\nu_{2,2}(a) = \Lambda_{W, b_2}(a_2)$.

The left hand side of \eqref{eq:sept-21-1} becomes
\begin{equation}
    \label{eq:sept-21-2}
    \mathbb{E}_{r \in [M_n]} (\nu(r) - 1) \mathbb{E}_{a \in [N'_n]^2} \prod_{j=1}^2 1_{A'_n}(a + r h_j) \theta_j(a + r h_j).
\end{equation}
As discussed in \cref{sec:tao-ziegler-summary}, the way we truncate $A_n$ (hence $A'_n$) allows us to shift $a + rh_1$ to $a$. Hence \eqref{eq:sept-21-2} is equal to 
\begin{equation}
    \mathbb{E}_{a \in [N'_n]^2} 1_{A'_n}(a) \theta_1(a) \mathbb{E}_{r \in [M_n]} (\nu(r) - 1) 1_{A'_n}(a + (h_2 - h_1) r) \theta_2(a + (h_2 - h_1) r).
\end{equation}
By Cauchy-Schwarz inequality, the square of previous expression is at most
\begin{multline}
\label{eq:sept-21-3}
    \mathbb{E}_{a \in [N'_n]^2} 1_{A'_n}(a) \theta_1(a) \times \\
    \times \mathbb{E}_{a \in [N_n]62} 1_{A'_n}(a) \theta_1(a) \left| \mathbb{E}_{r \in [M_n]} (\nu(r) - 1) 1_{A'_n}(a + (h_2 - h_1) r) \theta_2(a + (h_2 - h_1) r) \right|^2
\end{multline}
Observe that 
\[
    \mathbb{E}_{a \in [N'_n]^2} 1_{A'_n}(a) \theta_1(a) \leq \mathbb{E}_{a \in [N'_n]^2} \theta_1(a) = 1 + o(1)
\] 
where the last equality follows from the linear forms condition. Hence up to a factor of $1+o(1)$, \eqref{eq:sept-21-3} is at most
\begin{equation}
    \mathbb{E}_{a \in [N'_n]^2} \theta_1(a) \left| \mathbb{E}_{r \in [M_n]} (\nu(r) - 1) 1_{A'_n}(a + (h_2 - h_1) r) \theta_2(a + (h_2 - h_1) r) \right|^2
\end{equation}
By van der Corput's lemma, the previous expression is at most (up to a multiplicative constant)
\begin{multline}
    \mathbb{E}_{a \in [N'_n]^2} \theta_1(a) \mathbb{E}_{r \in [M_n], |s_1| < M_n} (\nu(r) - 1) (\nu(r + s_1) - 1) \times \\
    \times 1_{A'_n}(a + (h_2 - h_1) r) \theta_2(a + (h_2 - h_1) r) \times \\
    \times 1_{A'_n}( a + (h_2 - h_1) r + (h_2 - h_1) s_1) \theta_2(a + (h_2 - h_1) r + (h_2 - h_1) s_1).
\end{multline}
By shifting $a + (h_2 - h_1)r$ to $a$, the previous expression is equal to
\begin{multline}
\label{eq:sept-21-4}
    \mathbb{E}_{a \in [N'_n]^2} 1_{A'_n}(a) \theta_2(a) 
    1_{A'_n}( a + (h_2 - h_1) s_1) \theta_2(a +(h_2 - h_1) s_1)\\
    \mathbb{E}_{r \in [M_n], |s_1| < M_n} (\nu(r) - 1) (\nu(r + s_1) - 1) \theta_1(a + (h_1 - h_2) r).
\end{multline}
We now move everything that does not depend on $r$ outside of the average on $r$. Extra caution should be made here. At first glance, we may leave $\theta_1(a + (h_1 - h_2) r)$ inside since it seemingly  depends on $r$. However, as $h_1 - h_2 = (-1, 0)$
\begin{equation}
    \theta_1(a + (h_1 - h_2) r) = \nu_{1,1}(a_1 - r) \nu_{1,2}(a_2).
\end{equation}
Because $\nu_{1,2}(a_2)$ does not depend on $r$, we move it outside of the average on $r$. Thus \eqref{eq:sept-21-4} is equal to
\begin{multline}
     \mathbb{E}_{a \in [N'_n]^2} 1_{A'_n}(a) \theta_2(a) 
    1_{A'_n}( a + (h_2 - h_1) s_1) \theta_2(a +(h_2 - h_1) s_1) \nu_{1,2}(a_2)\\
    \mathbb{E}_{r \in [M_n], |s_1| < M_n} (\nu(r) - 1) (\nu(r + s_1) - 1) \nu_{1,1}(a_1 - r).
\end{multline}
By Cauchy-Schwarz inequality and van der Corput's lemma again, the square of above equation is at most (up to a multiplicative constant)
\begin{multline}
    \label{eq:sept-21-5}
    \mathbb{E}_{a \in [N'_n]^2, r \in [M_n], |s_1| < M_n, |s_2| < |M_n|} (\nu(r) - 1) (\nu(r + s_1) - 1) (\nu(r + s_2) - 1) (\nu(r + s_1 + s_2) - 1) \\
    \times \nu_{2,1}(a_1) \nu_{2,1}(a_1 + s_1)
    \nu_{1,2}(a_2) \nu_{1,1}(a_1 - r) \nu_{1, 1}(a_1 - r - s_2).
\end{multline}
Expanding out the last expression, we see that it is equal to a sum and difference of 16 averages of the form \eqref{eq:linearforms}.
Since no two forms in \eqref{eq:sept-21-5} are linearly dependent, each of these 16 averages is equal to $1+o(1)$ by the linear forms condition. Therefore $\eqref{eq:sept-21-5}$ is equal to $o(1)$. This finishes the proof of the special case.     

\subsection{The general case}
We now prove \cref{prop:1} in its full generality. Define $\nu = \Lambda_{W, 1}$ as before. For $1 \leq j \leq k$ and $1 \leq i \leq d$, define $\nu_{j,i}: \mathbb{Z}^d \to \mathbb{R}$ by
\begin{equation}
    \nu_{j,i}(a) = \begin{cases} \Lambda_{W, b_i}(a_i) \,\, \mbox{ if } h_{j,i} \neq h_{l, i} \; \forall 1 \leq l < j \\
    1 \mbox{ otherwise.}
    \end{cases}
\end{equation}
Then define $\theta_j = \prod_{i=1}^d \nu_{j,i}$. 

We need to show
\begin{equation}
    \mathbb{E}_{r \in [M_n]} (\nu(r) - 1) \mathbb{E}_{a \in [N'_n]^d} \prod_{j = 1}^k 1_{A'_n}(a + r h_j) \theta_j(a + r h_j) = o(1).
\end{equation}
As in the toy example, we perform the following steps:
\begin{enumerate}
    \item Shift $a$
    \item Move everything that does not depend on $r$ outside of the average on $r$
    \item Apply Cauchy-Schwarz inequality
    \item Bound $1_{A'_n} \theta_j$ by $\theta_j$
    \item Apply van der Corput's lemma
\end{enumerate}

After $k$ iterations, we get the following expression: 
\begin{multline}
    \label{eq:sept-21-6}
    \mathbb{E}_{a \in [N'_n]^d, r \in [M_n], |s_1| < M_n, \ldots, |s_k| < M_n} \prod_{R \subseteq [k]} \left( \nu(r + \sum_{l \in R} s_l) - 1 \right) \times \\
    \times \prod_{j=1}^k \prod_{i=1}^d \prod_{R \subseteq R_{j,i}} \nu_{j,i} \left( a_i + (h_{j,i} - h_{k,i}) r + \sum_{l \in R} (h_{j,i} - h_{l,i}) s_l \right)
\end{multline}
where $R_{j,i} = \{1 \leq l \leq k: h_{l,i} \neq h_{j,i}\}$. Again, this expression is a sum and difference of $2^{2^k}$ averages of the form \eqref{eq:linearforms}. In order to invoke the linear forms condition, it suffices to verify that no two forms appearing in \eqref{eq:sept-21-6} are linearly dependent.

Firstly, those forms appearing in $\nu$ are independent from one another because the appearance of $s_1, s_2, \ldots, s_k$ in each form is corresponding to a subset $R$ of $[k]$. They are also independent from the forms in $\nu_{j,i}$ because $a_i$ appears in $\nu_{j,i}$, but not in $\nu$.

Secondly, for $i_1 \neq i_2$ and any $j_1, j_2$ (not necessarily distinct), the forms appearing in $\nu_{j_1, i_1}$ and $\nu_{j_2, i_2}$ are independent because $a_{i_1}$ appears in $\nu_{j_1, i_1}$ while $a_{i_2}$ appears in $\nu_{j_2, i_2}$. 

For a fixed $i$, if $(h_{l,i} - h_{k,i}) r = (h_{j, i} - h_{k,i})r$ for some $l < j$, then $h_{l,i} = h_{j,i}$. By the way we define $\nu_{j,i}$, this would force $\nu_{j,i} = 1$. Hence if $\nu_{j,i} \neq 1$, $(h_{l,i} - h_{k,i}) r \neq (h_{j, i} - h_{k,i})r$. This implies the forms appearing in $\nu_{j,i}$ are independent from the forms appearing in $\nu_{l, i}$ for all $l < j$.

And lastly, for fixed $i$ and $j$, the forms appearing in $\nu_{j,i}$ are independent from one another because each form is in one-to-one correspondence with a subset of $R_{j,i}$. This finishes our proof.

\bibliographystyle{plain}
\bibliography{refs}

\end{document}